\documentclass[11pt]{article}
\usepackage{amssymb,amsmath,bm,graphicx}
\usepackage{mathrsfs}
\usepackage{natbib}
\usepackage{constants}
\usepackage{enumitem}
\usepackage{float}
\usepackage{hyperref}
\usepackage{atbegshi}% http://ctan.org/pkg/atbegshi
\AtBeginDocument{\AtBeginShipoutNext{\AtBeginShipoutDiscard}}
\usepackage{color}
\date{}
\begin{document}
\newcommand{\bea}{\begin{eqnarray}}
\newcommand{\ena}{\end{eqnarray}}
\newcommand{\beas}{\begin{eqnarray*}}
\newcommand{\enas}{\end{eqnarray*}}
\newcommand{\beq}{\begin{equation}}
\newcommand{\enq}{\end{equation}}
\def\qed{\hfill \mbox{\rule{0.5em}{0.5em}}}
\newcommand{\bbox}{\hfill $\Box$}
\newcommand{\ignore}[1]{}
\newcommand{\ignorex}[1]{#1}
\newcommand{\wtilde}[1]{\widetilde{#1}}
\newcommand{\qmq}[1]{\quad\mbox{#1}\quad}
\newcommand{\qm}[1]{\quad\mbox{#1}}
\newcommand{\nn}{\nonumber}
\newcommand{\Bvert}{\left\vert\vphantom{\frac{1}{1}}\right.}
\newcommand{\To}{\rightarrow}
\newcommand{\E}{\mathbb{E}}
\newcommand{\Var}{\mathrm{Var}}
\newcommand{\supp}{\mathrm{supp}}
\newcommand{\Cov}{\mathrm{Cov}}
\newcommand{\Corr}{\mathrm{Corr}}
\newcommand{\esssup}{\mathrm{ess} \ \mathrm{sup}}
\makeatletter
\newsavebox\myboxA
\newsavebox\myboxB
\newlength\mylenA
\newcommand*\xoverline[2][0.70]{%
    \sbox{\myboxA}{$\m@th#2$}%
    \setbox\myboxB\null% Phantom box
    \ht\myboxB=\ht\myboxA%
    \dp\myboxB=\dp\myboxA%
    \wd\myboxB=#1\wd\myboxA% Scale phantom
    \sbox\myboxB{$\m@th\overline{\copy\myboxB}$}%  Overlined phantom
    \setlength\mylenA{\the\wd\myboxA}%   calc width diff
    \addtolength\mylenA{-\the\wd\myboxB}%
    \ifdim\wd\myboxB<\wd\myboxA%
       \rlap{\hskip 0.5\mylenA\usebox\myboxB}{\usebox\myboxA}%
    \else
        \hskip -0.5\mylenA\rlap{\usebox\myboxA}{\hskip 0.5\mylenA\usebox\myboxB}%
    \fi}
\makeatother

\newtheorem{theorem}{Theorem}[section]
\newtheorem{corollary}[theorem]{Corollary}
\newtheorem{conjecture}[theorem]{Conjecture}
\newtheorem{proposition}[theorem]{Proposition}
\newtheorem{lemma}[theorem]{Lemma}
\newtheorem{definition}[theorem]{Definition}
\newtheorem{example}[theorem]{Example}
\newtheorem{remark}[theorem]{Remark}
\newtheorem{case}{Case}[section]
\newtheorem{condition}{Condition}[section]
\newcommand{\proof}{\noindent {\it Proof:} }

\title{{\bf\Large Concentration inequalities using approximate zero bias couplings with applications to Hoeffding's statistic under the Ewens distribution}}
\author{Nathakhun Wiroonsri \bigskip \\ Mathematics and Statistics with Applications Research Group \\
Department of Mathematics, King Mongkut's University of Technology Thonburi}
\thanks{This research was partially supported by NSTDA young researcher fund 64000343. \\
Email: nathakhun.wir@kmutt.ac.th}
\maketitle
\date{}

\begin{abstract} 
We prove concentration inequalities of the form $P(Y \ge t) \le \exp(-B(t))$ for a random variable $Y$ with mean zero and variance $\sigma^2$ using a coupling technique from Stein's method that is so-called approximate zero bias couplings. Applications to the Hoeffding's statistic where the random permutation has the Ewens distribution with parameter $\theta>0$ are also presented. A few simulation experiments are then provided to visualize the tail probability of the Hoeffding's statistic and our bounds. Based on the simulation results, our bounds work well especially when $\theta \le 1$.

\bigskip


\end{abstract}

\section{Introduction} \label{Intro}

Stein's method was first introduced by Charles Stein in his seminal paper (\cite{Stein72}) and is best known for obtaining non-asymptotic bounds on various distances for distributional approximations, see the text (\cite{CGS11}) and the introductory notes (\cite{Ross11}). Thus far, many applications in several areas such as combinatorics, statistics, statistical physics and applied sciences have been developed using this method. Meanwhile, some researchers have used Stein's method for non-distributional approximation contexts. \cite{Cha07} introduced a version of Stein's method for deriving concentration of measure inequalities. 

The name `concentration of measure' was first used in \cite{Tal95} when the author aimed to derive deviation bounds of the form $P(|X-\E X|\ge x)$ in the situations where the exact evaluation is impossible or is hard to obtain. A bound of this form is considered `good' if it decays sufficiently rapidly in $x$. For instance, most of the works in the literature seek for bounds that are decreasing at exponential rate.

Now we only restrict our attention to concentration inequalities using Stein's method. \cite{LN06} and \cite{Rat13} obtained non-uniform concentration inequalities for randomized orthogonal array sampling designs and random permutation sum, respectively, however, they are not in the form that we focus on. \cite{Cha07} used exchangeable pairs technique to provide concentration of measure for functions of dependent random variables, including the Hoeffding's statistic when the random permutation has the uniform distribution and the net magnetization in the Curie-Weiss model. Once the connection between Stein's method and concentration inequalities was uncovered, many results have been obtained using various techniques related to Stein's method. Concentration bounds using bounded size bias couplings was introduced in \cite{GG11a} and \cite{GG11b} with examples including the number of relatively ordered subsequences of a random permutation, sliding window statistics, the lightbulb process, and etc. The same type of results with the help of bounded zero bias couplings occurred in \cite{GI14} with applications to the Hoeffding's statistic where the random permutation has either the uniform distribution or one which is constant over permutations with the same cycle type and having no fixed points. The bounded size bias technique was strengthened in \cite{AB15} and the bounded assumption was relaxed by \cite{CGJ18} with applications to Erd\"{o}s-R\'{e}nyi graphs and permutation models.

In this work, we generalize the results in \cite{GI14} with zero bias couplings replaced by approximate zero bias couplings, first appeared in \cite{Wir17}. This coupling is constructed by adding a remainder term to the construction of the zero bias coupling generated by a $\lambda$-Stein pair (see \cite{GR97} or Section 4.4 of \cite{CGS11}). Here we provide brief detail, presented in \cite{Wir17}, of the construction of approximate zero bias couplings generated by approximate $\lambda,R$-Stein pairs. Recall that we call a pair of random variables $(Y',Y'')$, an \textit{approximate $\lambda,R$-Stein pair} if it is exchangeable and satisfies
\bea \label{expectvar}
\E Y' = 0, \quad \Var(Y') = \sigma^2 
\ena
with $\sigma^2 \in (0,\infty)$, and
\bea \label{approxsteinpair}
\E(Y''|Y') = (1-\lambda)Y'+R,
\ena
for some  $0< \lambda < 1$ and $R=R(Y')$. 

%\ncolor{
%Let $(Y',Y'')$ be an approximate $\lambda,R$-Stein pair.
%Taking expectation in \eqref{approxsteinpair}, using exchangeability and that $Y'$ has mean zero yields
%\beas
%0=\E Y''=(1-\lambda)\E Y'+\E R=\E R.
%\enas
%In addition, for any function $f$ such that the following expectations exist,
%\bea \label{approxsteinpair2}
%\E Y''f(Y') = \E(\E(Y''f(Y')|Y')) = \E(f(Y')\E(Y''|Y'))  = \E(f(Y')((1-\lambda)Y'+R)) \nn \\
%= (1-\lambda)\E Y'f(Y') + \E Rf(Y').
%\ena
%In particular, specializing \eqref{approxsteinpair2} to the case $f(y) = y$ yields
%\beas
%\E Y''Y' = (1-\lambda)\E Y'^2 + \E Y'R,
%\enas
%and thus
%\bea \label{diff2}
%\E(Y''-Y')^2 = 2\left(\E Y'^2 - \E Y''Y' \right) = 2\left(\lambda\E Y'^2 - \E Y'R \right) = 2\left(\lambda \sigma^2 -\E Y'R \right).
%\ena}

Now we state the following lemma proved in \cite{Wir17} which is the generalized version of the result in \cite{GR97} with a $\lambda$-Stein pair replaced by an approximate $\lambda,R$-Stein pair. We call $Y^*$ that satisfies \eqref{approxzerobias}, an \textit{approximate $Y'$-zero bias variable}.

\begin{lemma} \label{Wdagger}
Let $(Y',Y'')$ be an approximate $\lambda,R$-Stein pair with distribution $F(y',y'')$. Then when $(Y^\dagger,Y^\ddagger)$ has distribution 
\beas
dF^\dagger(y',y'') = \frac{(y''-y')^2}{\E (Y''-Y')^2}dF(y',y''),
\enas 
and $U$ has a uniform distribution on $[0,1]$, independent of $Y^\dagger,Y^\ddagger$, the variable $Y^* = UY^\dagger + (1-U)Y^\ddagger$ satisfies
\bea \label{approxzerobias}
\E [Y'f(Y')] = \sigma^2\E f'(Y^*) - \frac{\E Y'R}{\lambda}\E f'(Y^*) + \frac{\E Rf(Y')}{\lambda}
\ena
for all absolutely continuous functions $f$.
\end{lemma}

The core idea of Stein's method is based on the fact that a mean zero random variable $W$ is normal if and only if $\E[Wf(W)]=\sigma^2\E[f'(W)]$ for all absolutely continuous function. Therefore if $Y'$ and $Y^*$ can be closely coupled by using Lemma \ref{Wdagger} satisfying \eqref{approxzerobias} with some small remainder term $R$ then we expect that the behavior of $Y'$, including the decay of its tail probabilities, may be simlar to that of the normal. \cite{Wir17} obtained $L^1$ and $L^\infty$ bounds between $Y'/\sigma$ and the standard normal in term of the difference $|Y'-Y^*|$ and the remainder term $R$. Here we focus on the concentration bounds on the same setting.

We also apply our result to the Hoeffding's statistic (\cite{Hoe51}) where the random permutation has the Ewens distribution. That is, we study the distribution of 
\bea \label{combdef}
Y =\sum_{i=1}^n a_{i,\pi(i)}
\ena
where $A \in \mathbb{R}^{n \times n}$ is a given real matrix with components $\{a_{i,j} \}_{i,j=1}^n$ and $\pi \in \mathcal{S}_n$ has the Ewens distribution where $\mathcal{S}_n$ denotes the symmetric group. We note that the distribution of $Y$ is of interest as it has a connection to nonparametric tests in statistics. \cite{Wir17} has succeeded in providing $L^1$ and $L^\infty$ bounds of order $\sigma^{-1}$ between $Y$ and the normal by constructing an approximate $Y$-zero bias coupling and deriving the bounds based on the coupling. In the present work, we follow the construction there and apply the following main results to obtain concentration inequalities.

We now state the main theorem and then add a remark comparing the results here to the ones in \cite{GI14}. 

\begin{theorem} \label{main1}
Let $Y$ be a mean zero random variable with variance $\sigma^2 >0$ and moment generating function $m(s) = \E[e^{sY}]$, and let $Y^*$ have the approximate $Y$-zero bias distribution satisfying \eqref{approxzerobias}, with $R=R(Y)$ and $0<\lambda<1$. Let $B_1=\E|R|/\lambda$ if $|R|$ and $e^{sY}$ are negatively correlated for $0<s<\alpha$ where $\alpha$ is given differently for each bound below and $B_1 = \esssup|\E[R|Y]|/\lambda$ otherwise and $B_2 = |\E(YR)|/\lambda$. Assume that $B_1$ and $B_2$ are finite. 
\begin{enumerate}[label=(\alph*)]
	\item If $Y^*-Y \le c$ for some $c>0$ and $m(s)$ exists for all $s \in[0,\alpha)$ with $\alpha=1/c$, then for all $t > 0$
\bea \label{mainbound1}
P(Y \ge t) \le \exp\left(-\frac{t(t-2B_1)}{2(\sigma^2+B_2+ct)}\right).
\ena
The same upper bound holds for $P(Y \le -t)$ if $Y-Y^* \le c$ when $m(s)$ exists for all $s \in (-\alpha,0]$. If $|Y^*-Y| \le c$ for some $c>0$ and $m(s)$ exists for all $s \in [0,\alpha)$ with $\alpha=2/c$ then for all $t >0$
\bea \label{mainbound2}
P(Y \ge t) \le \exp\left(-\frac{t(t-2B_1)}{10(\sigma^2+B_2)/3+ct}\right),
\ena
with the same upper bound for $P(Y\le -t)$ if $Y-Y^* \le c$ and $m(s)$ exists in $(-\alpha,0]$.
\item If $Y^*-Y \le c$ for some $c>0$ and $m(s)$ exists at $\alpha = (\log (t-B_1)-\log \log (t-B_1))/c$ then for $t>B_1$,
\bea \label{mainbound3}
P(Y \ge t) &\le& \exp\left(-\frac{t-B_1}{c}\left(\log (t-B_1)  - \log \log (t-B_1) -\frac{\sigma^2+B_2}{c}\right)\right) \nn\\
   &\le& \exp\left(-\frac{t-B_1}{2c}\left(\log (t-B_1) -\frac{2(\sigma^2+B_2)}{c}\right)\right).
\ena
If $Y-Y^*\le c$ then the same bound holds for the left tail $P(Y \le -t)$ when $m(-\alpha)<\infty$.
\end{enumerate}
\end{theorem}
Part (b) shows that the respective asymptotic orders as $t \rightarrow \infty$ of $\exp(-t/(2c))$ and $\exp(-t/c)$ of the bounds \eqref{mainbound1} and \eqref{mainbound2} can be improved to $\exp(-t\log t/(2c))$.

\begin{remark} \label{rem1}
\begin{enumerate}
	\item If $R=0$, then $B_1$ and $B_2$ are zero and therefore the bounds \eqref{mainbound1}, \eqref{mainbound2} and \eqref{mainbound3} in Theorem \ref{main1} are exactly the same as (2), (3) and (4) of \cite{GI14}, respectively. See Remark 2.1 of \cite{GI14} and a few paragraphs thereafter for a comparison to previous results in the literature.
	\item We note that the bounds in \eqref{mainbound1} and \eqref{mainbound2} will be effective only for $t > 2B_1$, and the bound \eqref{mainbound3} for $t>B_1$ where $B_1$ is given in the statement of the theorem. Therefore the three bounds in Theorem \ref{main1} will be useful only if we can construct an approximate $\lambda,R$-Stein pair in such a way that $B_1$ or $2B_1$ is in the range of the support of $Y$. Otherwise, the bounds will be trivial. 
\end{enumerate}

\end{remark}

Stein's method has been used with Ewens distribution recently in \cite{HR21} where they obtained an upper bound on the total variabtion distance between the random partition generated by Ewens Sampling Formula and the Poisson-Dirichlet distribution. Ewens distribution has also attracted attention in other fields of mathematics for a while (see \cite{JKB97}, \cite{ABT03}, \cite{Cra16} and \cite{SJT20} for instance.)

The remainder of this work is organized as follows. We present an application to the Hoeffding's statistic where the random permutation has the Ewens distribution in Section \ref{app}. A few simulation experiments are then performed in Section \ref{sim}. Finally, we devote  the last section for the proof of our main theorem.

\section{The Hoeffding's statistic under the Ewens measure} \label{app}

As introduced in the first section, we study the distribution of 
\bea \label{combdef}
Y =\sum_{i=1}^n a_{i,\pi(i)}
\ena
where $A \in \mathbb{R}^{n \times n}$ is a given real matrix with components $\{a_{i,j} \}_{i,j=1}^n$ and $\pi \in \mathcal{S}_n$ has the Ewens distribution. The original work (\cite{Hoe51}) studied the asymptotic normality of $Y$ when $\pi \in \mathcal{S}_n$ has a uniform distribution. We note that, in this section, we consider the symbols $\pi$ and $Y$ interchageable with $\pi'$ and $Y'$, respectively.

We first briefly describe the Ewens distribution and state some important properties. The Ewens distribution ${\cal E}_\theta$ on the symmetric group ${\cal S}_n$ with parameter $\theta > 0$, was first introduced in \cite{Ewe72} and used in population genetics to describe the probabilities associated with the number of times that different alleles are observed in the sample; see also \cite{ABT03} for the description in mathematical context. Let $\mathbb{N}_k=[k,\infty) \cap \mathbb{Z}$. In the following, for $x \in \mathbb{R}$ and $n \in \mathbb{N}_1$, we use the notations
\beas
x^{(n)} = x(x+1)\cdots(x+n-1)  \text{ \ \ and \ \ } x_{(n)} = x(x-1)\cdots(x-n+1).
\enas
For a permutation $\pi \in \mathcal{S}_n$, the Ewens measure is given by 
\bea \label{ewensdef}
P_{\theta}(\pi) = \frac{\theta^{\#(\pi)}}{\theta^{(n)}},
\ena
where $\#(\pi)$ denotes the number of cycles of $\pi$. We note that ${\cal E}_\theta$ specializes to the uniform distribution over all permutations when $\theta=1$.

A permutation $\pi_n \in \mathcal{S}_n$ with the distribution ${\cal E}_\theta$ can be constructed by the ``so called'' the Chinese restaurant process (see e.g. \cite{Ald85} and \cite{Pit96}), as follows. For $n=1$, $\pi_1$ is the unique permutation that maps $1$ to $1$ in $\mathcal{S}_1$. For $n \ge 2$, we construct $\pi_n$ from $\pi_{n-1}$ by either adding $n$ as a fixed point with probability $\theta/(\theta+n-1)$, or by inserting $n$ uniformly into one of $n-1$ locations inside a cycle of $\pi_{n-1}$, so each with probability $1/(\theta+n-1)$.
Following this construction, for $i,j,k,l \in [n]$ where $[n] = \{1,2,\ldots,n \}$ such that $i \ne j$ and $k\ne l$, we have
\bea \label{eq:marginal.1}
P_\theta(\pi(i) = k) =
\begin{cases}
\frac{1}{\theta+n-1} \text{ \ if \ } k \ne i,\\
\frac{\theta}{\theta+n-1} \text{ \ if \ } k = i
\end{cases}.
\ena

A distribution on $\mathcal{S}_n$ is said to be constant on cycle type if the probability of any permutation $\pi \in {\cal S}_n$ depends only on the cycle type $(c_1,\ldots,c_n)$ where $c_q(\pi)$ is the number of $q$ cycles of $\pi$ and we write $c_q$ for $c_q(\pi)$ for simplicity. Concentration inequalities of $Y$ where $\pi \in \mathcal{S}_n$ has either a uniform distribution or one which is constant over permutations with the same cycle type and having no fixed points were obtained in \cite{GI14} by using zero bias couplings. It follows from \eqref{ewensdef} directly that the Ewens distribution is constant on cycle type and allows fixed point. Due to the difference, the zero biasing technique does not work here and thus we use approximate zero biasing as discussed in the first section instead. We first follow the construction of approximate $Y$-zero bias couplings presented in \cite{Wir17} and then apply Theorem \ref{main1} to obtain concentration bounds. It is interesting to study the robustness and the sensitivity of the tail probabilities when the controlled assumptions have changed. 

We begin with stating the definitions, notations and assumptions used in \cite{Wir17}. Letting
\bea \label{adotdot}
a_{\bullet,\bullet} = \frac{1}{n(\theta+n-1)} \left(\theta \sum_{i=1}^n a_{i,i}+ \sum_{i \ne j}a_{ij}\right),
\ena 
by \eqref{eq:marginal.1}, we have
\bea \label{meanY}
\E Y = n a_{\bullet,\bullet}.
\ena
Letting
\beas
\widehat{a}_{i,j} = a_{i,j} - a_{\bullet,\bullet},
\enas
and using \eqref{meanY}, we have 
\beas
\E \left[\sum_{i=1}^n \widehat{a}_{i,\pi(i)}\right]=\E \left[\sum_{i=1}^n (a_{i,\pi(i)}-a_{\bullet,\bullet})\right] = \E \left[\sum_{i=1}^n a_{i,\pi(i)} - na_{\bullet,\bullet}\right] = 0.
\enas
As a consequence, replacing $a_{i,j}$ by $\widehat{a}_{i,j}$, we assume without loss of generality that
\bea \label{zeroassumption}
\E Y = 0 \text{ \ and \ } a_{\bullet,\bullet} = 0, 
\ena 
and for simplicity, we consider only the symmetric case, that is, $a_{i,j} = a_{j,i}$ for all $i,j$. Also, to rule out trivial cases, we assume in what follows that $\sigma^2>0$. The explicit form of $\sigma^2$ was calculated in Lemma 4.8 of \cite{Wir17} and as discussed in Remark 4.9 of \cite{Wir17}, $\sigma^2$ is of order $n$ when the elements of $A$ are well chosen so that they do not depend on $n$ and most of the elements are not the same number. In this case, there exists $N \in \mathbb{N}_1$ such that $\sigma^2>0$ for $n > N$.

The following theorem presents concentration inequalities for $Y$.

\begin{theorem} \label{appthm}
Let $n \ge 6$ and let $A = \{a_{i,j} \}_{i,j=1}^n$ be an array of real numbers satisfying
\beas
a_{i,j} = a_{j,i}.
\enas
Let $\pi \in \mathcal{S}_n$ be a random permutation with the distribution ${\cal E}_\theta$, with $\theta > 0$. 
Then, with $Y$ the sum in \eqref{combdef} and
\beas
M = \max_{i,j} |a_{i,j} -a_{\bullet,\bullet}|, \  a_{\bullet,\bullet} \text{ \ as in \eqref{adotdot}},
\enas
the bounds in Theorem \ref{main1} hold with $Y$ replaced by $Y-\E Y$ and $c = 20 M$. 

Moreover, 
\beas
B_1 \le (6n+4.8\theta) M  \text{ \ and \ } 
B_2 \le \left(3\kappa_{\theta,n,1} + 1.2 \theta + \frac{1.2(\kappa_{\theta,n,1}(\theta+1)+\kappa_{\theta,n,2})}{n}  \right)M \sigma,
\enas
where
\bea \label{kappadef}
\kappa_{\theta,n,1} = \sqrt{\frac{\theta^2n_{(2)}}{(\theta+n-1)_{(2)}} + \frac{\theta n}{\theta+n-1}}
\ena
and
\bea \label{kappadef2}
\kappa_{\theta,n,2} = \sqrt{\frac{\theta^4n_{(4)}}{(\theta+n-1)_{(4)}}+\frac{4\theta^3 n_{(3)}}{(\theta+n-1)_{(3)}}+\frac{2\theta^2n_{(2)}}{(\theta+n-1)_{(2)}}}.
\ena

In particular, if $A$ is chosen so that $|R|$ given in \eqref{Rdef} and and $e^{sY}$ are negatively correlated for $0<s<\alpha$, then 
\beas
B_1 \le \frac{\theta M (3.6n+2.4\theta -3)}{\theta+n-1} + \frac{ \theta^2 n M}{2(\theta+n-1)_{(2)}}.
\enas
\end{theorem}

As mentioned earlier in Remark \ref{rem1}, the bounds arised in Theorem \ref{appthm}, are effective for $t \ge 2B_1$ where $B_1$ depends on the remainder term of the approximate $\lambda,R$-Stein pair and the bounds are trivial outside this range. Hence if $2B_1$ is out of the support of $Y-\E Y$, then the bound will not be useful.  The best possible bound of $B_1$ that we can obtain in general is $(6n+4.8\theta) M$ which depends on $n$, nevertheless, it becomes $\theta M (3.6n+2.4\theta -3)/(\theta+n-1) +  \theta^2 n M/(2(\theta+n-1)_{(2)})$ which is of constant order in the case that $|R|$ and and $e^{sY}$ are negatively correlated for $0<s<\alpha$.

To prove Theorem \ref{appthm}, we apply Theorem \ref{main1} from the introduction. That is, we have to construct an approximate $Y$-zero bias coupling $Y^*$ from an approximate $\lambda,R$-Stein pair $(Y,Y'')$ by using Lemma \ref{Wdagger} in such a way that $|Y^*-Y|$ and $\E[R|Y]$ are bounded and that $|\E(YR)|<\infty$. For this purpose, we state the following two lemmas. The first, Lemma \ref{combapproxstein} proved in \cite{Wir17}, constructs an approximate $4/n,R$-Stein pair $(Y',Y'')$ for $n\ge 6$ where $Y'$ is given by \eqref{combdef} with $\pi'$ replacing $\pi$. The bounds related to the remainder term $R$ are then obtained in Lemma \ref{remainderbounds}. Below, for $i,j \in [n]$, we write $i\sim j$ if $i$ and $j$ are in the same cycle, let $|i|$ be the length of the cycle containing $i$ and let $\tau_{i,j}$, $i,j \in [n]$ be the permutation that transposes $i$ and $j$.

\begin{lemma}[\cite{Wir17}] \label{combapproxstein}
For $n \ge 6$, let $\{a_{i,j}\}_{i,j=1}^n$ be an array of real numbers satisfying $a_{i,j} = a_{j,i}$ and $a_{\bullet,\bullet}=0$ where $a_{\bullet,\bullet}$ is as in \eqref{adotdot}. Let $\pi \in \mathcal{S}_n$ a random permutation has the Ewens measure ${\cal E}_\theta$ with $\theta>0$, and let $Y'$ be given by \eqref{combdef}. Further, let $I,J$ be chosen independently of $\pi$, uniformly from all pairs of distinct elements of $\{1,\ldots n\}$. Then, letting $\pi'' = \tau_{I,J}\pi \tau_{I,J}$ and $Y'$ be given by \eqref{combdef} with $\pi'$ replacing $\pi$, $(Y',Y'')$ is an approximate $4/n,R$-Stein pair with 
\bea \label{Rdef}
R(Y') =  \frac{1}{n(n-1)} \E [T|Y']
\ena
where
\bea \label{eq:defT}
T &=& 2(n+c_1-2(\theta+1)) \sum_{|i|=1}a_{i,i} - 2(c_1-2\theta)\sum_{|i| \ge 2} a_{i,i} \nn\\
 && \hspace{80pt}-4\sum_{|i|=1,|j|=1, j \ne i}a_{i,j} -4\sum_{|i|=1,|j|\ge 2}a_{i,j}.
\ena
\end{lemma}

\begin{lemma} \label{remainderbounds}
Let $(Y',Y'')$ be an approximate $4/n,R$-Stein pair constructed as in Lemma \ref{combapproxstein} with $R$ as in \eqref{Rdef}. Then
\bea \label{Rbound}
|\E[R|Y]| \le \frac{10n+8\theta}{n-1} M  \text{ \ \ a.s},
\ena 
\bea \label{Rbound2}
\E|R| \le \frac{\theta M (12n+8\theta -10)}{(n-1)(\theta+n-1)} + \frac{2 \theta^2 M}{(\theta+n-1)_{(2)}},
\ena 
and
\bea \label{YRbound}
|\E Y'R| \le \left(10\kappa_{\theta,n,1} + 4 \theta + \frac{4(\kappa_{\theta,n,1}(\theta+1)+\kappa_{\theta,n,2})}{n}  \right) \frac{ M \sigma}{n-1}
\ena
where $M = \max_{i,j} |a_{i,j} -a_{\bullet,\bullet}|$ with $a_{\bullet,\bullet}$ as in \eqref{adotdot} and $\kappa_{\theta,n,1}$ and $\kappa_{\theta,n,2}$ are given in \eqref{kappadef} and \eqref{kappadef2}, respectively.
\end{lemma}
\proof
The bounds \eqref{Rbound2} of $\E|R|$ and \eqref{YRbound} of $|\E Y'R|$ were shown in \cite{Wir17}. To prove \eqref{Rbound}, since $|\E[R|Y]| = |\E[T|Y]|/(n(n-1)) \le |T|/(n(n-1))$, it suffices to bound $T$ in \eqref{eq:defT}. For any fixed $\pi \in \mathcal{S}_n$, using that $|\{i\in[n]:|i|=1 \}| = c_1$ and $|\{i\in[n]:|i|\ge 2 \}| = n-c_1$  we have
\beas
|T| \le 4\theta n M \text{ \ \  a.s. \ \  if \ \ } c_1=0,
\enas
and, if $c_1 > 0$,
\beas
|T| \le (2n+2c_1+4\theta+4)c_1M + (2c_1+4\theta)(n-c_1)M \\ 
+ 4c_1(c_1-1)M + 4c_1(n-c_1)M \text{ \ \ a.s. }
\enas  
Using that $c_1 \le n$ and that $n(n-c_1)$ is maximized when $c_1 = n/2$, it follows that
\beas
|T| \le (10n^2+8\theta n) M  \text{ \ \ a.s.}
\enas
\bbox

Now we have all ingredients to prove Theorem \ref{appthm}.
\bigskip

\noindent {\bf Proof of Theorem \ref{appthm}}:
First we construct an approximate $4/n,R$-Stein pair $(Y',Y'')$ as in Lemma \ref{combapproxstein} with the remainder $R$ given in \eqref{Rdef}. Then, constructing an approximate $Y$-zero bias variable $Y^*$ from the approximate Stein pair $(Y',Y'')$ as in the proof of Theorem 4.1 of \cite{Wir17}, we have $|Y^*-Y| \le 20M$. Invoking Theorem \ref{main1} with $\lambda = 4/n$ using the bounds of $|R|$ and $|\E YR|$ in \eqref{Rbound} and \eqref{YRbound} of Lemma \ref{remainderbounds}, respectively, completes the proof.
\bbox 

\bigskip

\section{Simulation} \label{sim}

In this section, we perform 4 simulation experiments. $A$ is generated by first randomly choosing $X$ from $N(1,0.2)$ and $N(-1,0.2)$ equally and then letting $B = X+X^T$ and $A = B-b_{\bullet,\bullet}$ where $b_{\bullet,\bullet}$ is defined similarly to \eqref{adotdot}. We intentionally select $X$ such that $\Cov(|R|,e^{sY})<0$ for all $0<s<\alpha$ with $\alpha$ as defined in Theorem \ref{main1} so that the negative correlation assumption there is satisfied. To generate Ewens permutations with $\theta \ne 1$, we apply the accept-reject algorithm (see \cite{RC04}). Recall that the discrete version of the accept-reject algorithm is stated as follow.

\begin{theorem} [\cite{Neu51}]
Let $Y$ and $V$ be discrete random variables with probability mass functions $f_Y(y)$ and $f_V(v)$, respectively, where $f_Y$ and $f_V$ have common support. Define
\beas
C = \sup_y \frac{f_Y(y)}{f_V(y)} < \infty.
\enas 
To generate a random variable $Y \sim f_Y$:
\begin{enumerate}[label=(\alph*)]
	\item Generate $U\sim$ Uniform$(0,1)$, $V\sim f_V$, independently.
	\item If $U \le \frac{f_Y(V)}{Cf_V(V)}$, set $Y=V$; otherwise, return to step (a).
\end{enumerate}
\end{theorem}

In our case, $Y$ and $V$ are $\mathcal{E}_\theta$ with $\theta\ne 1$ and uniform, respectively, with common support be all possible permutations $\pi \in \mathcal{S}_n$. Recall that the probability mass functions of $Y$ and $V$ are
\beas
f_Y(\pi) = \frac{\theta^{\#(\pi)}}{\theta^{(n)}} \text{ \ and \ } f_V(\pi) = \frac{1}{n!},
\enas
where $\#(\pi)$ denotes the number of cycles of $\pi$. It is easy to see that
\beas
C = 
\begin{cases}
\frac{n!\theta}{\theta^{(n)}} \text{ \ if \ } \theta<1 \\
\frac{n!\theta^n}{\theta^{(n)}} \text{ \ if \ } \theta>1 \\
\end{cases}.
\enas

\subsection{Experiment 1}
\begin{figure}[H]
\centering
\includegraphics[width=3.5in]{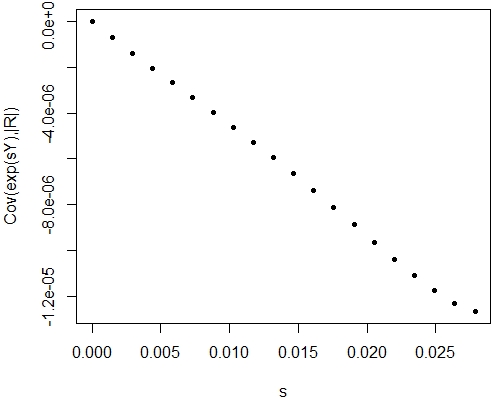}
\caption{$Cov(e^{sY},|R|)$ for $0\le s\le 2/c$.}
\label{cov0}
\end{figure}

\begin{figure}[H]
\centering
\includegraphics[width=3.5in]{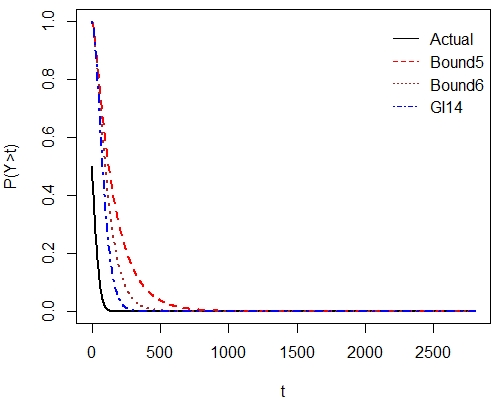}
\caption{Plot of $P(Y>t)$ when $n=1000$, $\theta=1$.}
\label{ex0}
\end{figure}

Let $n=1000$ and $\theta=1$. By simulating $Y$ and $R$ for 1,000,000 times, we have $\sigma^2=2082.80$, $B_1=0.74$, $B_2=2.13$, $M=3.59$, support$=[-2795.86,2803.19]$ and $\Cov(Y,|R|)=-0.0005$ from the sample. Figures \ref{cov0} and \ref{ex0} show $Cov(e^{sY},|R|)$ for $0<s<\alpha$ and the plot of $P(Y>t)$ compared with the bounds from \eqref{mainbound1} and \eqref{mainbound2}, respectively. Note that we intend not to show \eqref{mainbound3} since $Cov(e^{sY},|R|)>0$ at $\alpha = (\log (t-B_1)-\log \log (t-B_1))/c$ for large $t$. Since $\theta=1$ is equivalent to the uniform distribution on $\pi \in \mathcal{S}_n$, we also plot the bound from \cite{GI14}. It is clear that the bound from \cite{GI14} is smaller as our coupling construction yields larger $|Y-Y^*|$ and that our bound is derived in a way to handle the remainder term so it is not as tight. However, the bound from \cite{GI14} does not work with $\theta \ne 1$. 

\subsection{Experiment 2}

\begin{figure}[H]
\centering
\includegraphics[width=3.5in]{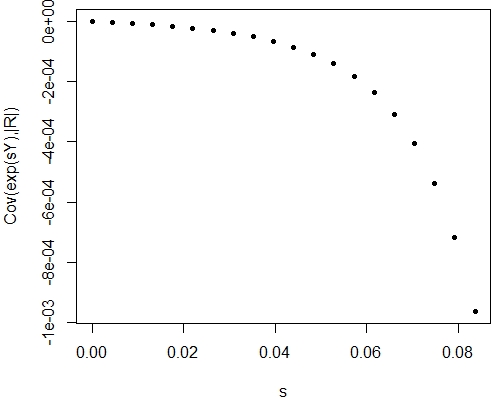}
\caption{$Cov(e^{sY},|R|)$ for $0\le s\le 2/c$.}
\label{cov2}
\end{figure}

\begin{figure}[H]
\centering
\includegraphics[width=3.5in]{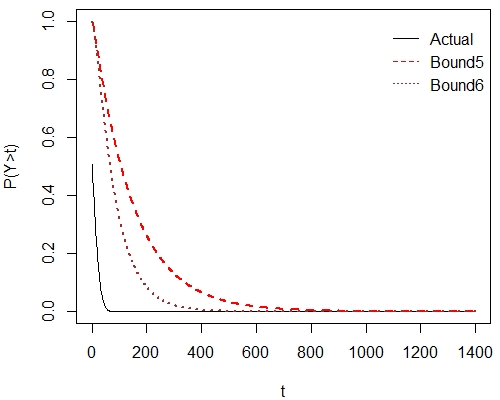}
\caption{Plot of $P(Y>t)$ when $n=1000$, $\theta=0.8$.}
\label{ex2}
\end{figure}

Since generating Ewens permutations with $\theta \ne 1$ requires an algorithm, we are not able to generate a very big sample with large $n$. Let $n=1000$ and $\theta=0.8$. $\theta<1$ implies that the chance of having a big number of cycles is less likely. By simulating $Y$ and $R$ 10,000 times, we have $\sigma^2=507.20$, $B_1=0.29$, $B_2=0.28$, $M=1.60$, support$=[-1400.45,1398.01]$ and $\Cov(Y,|R|)=-0.001$ from the sample. Figures \ref{cov2} and \ref{ex2} show $Cov(e^{sY},|R|)$ for $0<s<\alpha$  and the plot of $P(Y>t)$ compared with the bounds from \eqref{mainbound1} and \eqref{mainbound2}, respectively. Note that the bound \eqref{mainbound3} is too large in this case because of the term $\frac{2(\sigma^2+B_2)}{c}$.

\subsection{Experiment 3}
\begin{figure}[H]
\centering
\includegraphics[width=3.5in]{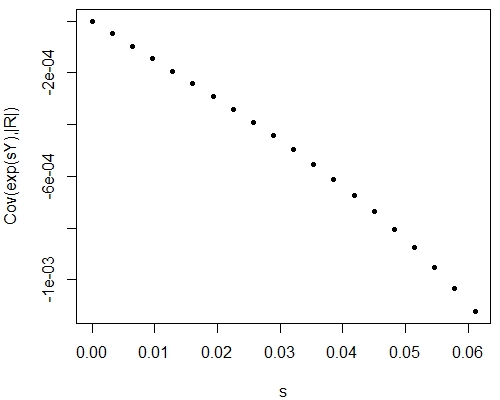}
\caption{$Cov(e^{sY},|R|)$ for $s = (\log (t-B_1)-\log \log (t-B_1))/c$.}
\label{cov3}
\end{figure}

\begin{figure}[H]
\centering
\includegraphics[width=3.5in]{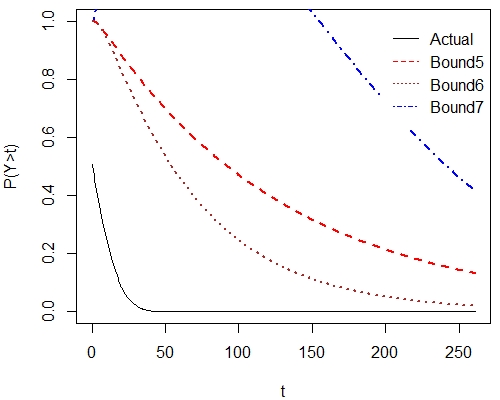}
\caption{Plot of $P(Y>t)$ when $n=100$, $\theta=1.05$.}
\label{ex3}
\end{figure}

When increasing $\theta$, it is more expensive to generate a permutation. Now we let $n=100$ and $\theta=1.05$. $\theta> 1$ provides that the chance of having a big number of cycles is more likely. By simulating $Y$ and $R$, 10,000 times, we have $\sigma^2=212.75$, $B_1=0.88$, $B_2=2.42$, $M=3.15$, support$=[-255.59,260.86]$ and $\Cov(Y,|R|)=-0.015$ from the sample. Figures \ref{cov3} and \ref{ex3} show $Cov(e^{sY},|R|)$ for $0<s<\alpha$  and the plot of $P(Y>t)$ compared with the bounds from \eqref{mainbound1}, \eqref{mainbound2} and \eqref{mainbound3}, respectively.

\subsection{Experiment 4}
\begin{figure}[H]
\centering
\includegraphics[width=3.5in]{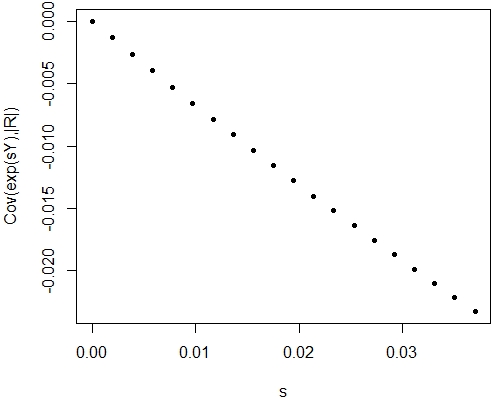}
\caption{$Cov(e^{sY},|R|)$ for $s = (\log (t-B_1)-\log \log (t-B_1))/c$.}
\label{cov4}
\end{figure}

\begin{figure}[H]
\centering
\includegraphics[width=3.5in]{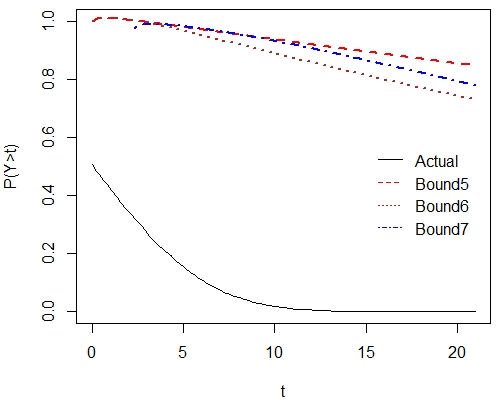}
\caption{Plot of $P(Y>t)$ when $n=10$, $\theta=2$.}
\label{ex4}
\end{figure}

For this experiment, we intend to illustrate the case that $\theta$ is much greater than 1 which results in a permutation with much larger number of cycles. However, the simulation is very expensive which even takes 93.41 iterations on average to get only one permutation with $n=10$ and $\theta=2$. Based on our simulation result with sample size of 10,000, the mean of the number of cycles is 4.03. By simulating $Y$ and $R$ 10,000 times, we have $\sigma^2=23.14$, $B_1=1.27$, $B_2=4.60$, $M=2.70$, support$=[-20.79,20.95]$ and $\Cov(Y,|R|)=-0.69$ from the sample. Figures \ref{cov4} and \ref{ex4} show $Cov(e^{sY},|R|)$ for $0<s<\alpha$ and the plot of $P(Y>t)$ compared with the bounds from \eqref{mainbound1}, \eqref{mainbound2} and \eqref{mainbound3}, respectively.

\subsection{Conclusion from the simulation results}

We can notice from the four experiments that our bounds seem to be more useful for smaller $\theta$, especially when $\theta \le 1$. For $\theta > 1$, the bounds are not sufficiently small even at the upper limit of the support. However, our experiments with $\theta>1$ are only done with small $n$. Therefore, it is interesting to see whether our bounds will work better with larger $n$.

Also based on our experiments, the bound \eqref{mainbound2} dominates the other in all cases. Since the bound \eqref{mainbound3} actually has the fastest decaying rate at $\exp\left(-t\log t\right)$, the bound \eqref{mainbound3} should be better if we consider $Y$ with a larger support.

\section{Proof of the main result} \label{proof}
In this section, we follow the same technique as in \cite{GI14} but handle a few extra terms that include $R$. 

\bigskip

\noindent {\bf Proof of Theorem \ref{main1}}:
Let $m(s) = \E[e^{sY}]$ and $m^*(s) = \E[e^{sY^*}]$. Using that $Y^*-Y \le c$, for $s \ge 0$ we have
\bea \label{mrelate}
m^*(s) = \E[e^{s(Y^*-Y)}e^{sY}] \le \E[e^{cs}e^{sY}] = e^{cs}m(s).
\ena

(a) If $m(s)$ exists in an open interval containing $s$ we may interchange expectation and differentiation at $s$ to obtain
\bea \label{m'bound}
m'(s) &=& \E[Ye^{sY}] = \sigma^2\E[se^{sY^*}] -\frac{\E YR}{\lambda}\E[se^{sY^*}] + \frac{\E [ R e^{sY} ]}{\lambda} \nn\\
      &=& \sigma^2 s m^*(s) - \frac{\E YR}{\lambda} s m^*(s) +  \frac{\E [ R e^{sY} ]}{\lambda},
\ena
where we have used the approximate zero bias relation in the second equality.

For $\alpha \in (0,1/c)$ and $0 \le s \le \alpha$, using $e^x \le 1/(1-x)$ for $x \in [0,1)$ we have
\bea \label{expbound}
m^*(s) \le e^{cs} m(s) \le \frac{1}{1-sc}m(s).
\ena
Using \eqref{m'bound} and applying \eqref{expbound}, we have
\bea \label{m'mineq}
m'(s) &\le& \frac{\sigma^2s}{1-sc}m(s) + \frac{|\E(YR)|}{\lambda} \frac{s}{1-sc} m(s) + \frac{|\E(R e^{sY})|}{\lambda} \nn\\
      &=& \frac{\sigma^2s}{1-sc}m(s) + \frac{|\E(YR)|}{\lambda} \frac{s}{1-sc} m(s) + \frac{|\E[e^{sY}\E[R|Y] ]|}{\lambda} \nn\\
			&\le&  \frac{(\sigma^2+|\E(YR)|/\lambda)s}{1-sc}m(s)+ \frac{\E[e^{sY}|\E[R|Y]|]}{\lambda}.
\ena
Recalling from the statement of the theorem that $B_1=\E|R|/\lambda$ if $|R|$ and and $e^{sY}$ are negatively correlated for $0<s<\alpha$  and $B_1 = \esssup|\E[R|Y]|/\lambda$ otherwise and $B_2 = |\E(YR)|/\lambda$, it follows from the last expression that 
\beas
m'(s) \le \frac{(\sigma^2+B_2)s}{1-sc}m(s)+ B_1 m(s).
\enas
Dividing both sides by $m(s)$, integrating over $[0,\alpha]$ and using that $m(0) = 1$ yield
\beas
\log(m(\alpha)) = \int_0^\alpha \frac{m'(s)}{m(s)} ds 
             \le \frac{\sigma^2+B_2}{1-\alpha c}\int_0^\alpha s ds
                                     + B_1 \int_0^\alpha 1 ds 
					   = \frac{(\sigma^2+B_2) \alpha^2}{2(1-\alpha c)} + B_1\alpha,
\enas
and therefore
\beas
m(\alpha) \le \exp\left(\frac{(\sigma^2+B_2)\alpha^2}{2(1-\alpha c)}  + B_1\alpha\right).
\enas

We note that the claim for $t=0$ is clear as the right hand side of \eqref{mainbound1} is one. Applying Markov's inequality for $t>0$, we obtain
\beas
P(Y \ge t) = P(e^{\alpha Y} \ge e^{\alpha t}) \le e^{-\alpha t} m(\alpha) \le \exp\left(-\alpha (t-B_1) + \frac{(\sigma^2+B_2) \alpha^2 }{2(1-\alpha c)}\right).
\enas
Letting $\alpha = t/(\sigma^2+B_2+ct)$ that lies in $(0,1/c)$, \eqref{mainbound1} holds. The result for the left tail when $Y-Y^* \le c$ follows immediately from the same proof by letting $X=-Y$ and the fact that $X^* = -Y^*$ has the approximate $X$-zero bias distribution with $R$ replaced by $-R$. We also note that the results for the left tail in the other cases below follow by the same argument so it suffices to prove only the results for the right tail.

Next we prove \eqref{mainbound2}. Following the same argument as in \cite{GI14} and using that
\beas
e^y-e^x \le \frac{|y-x|(e^y+e^x)}{2} \text{ \ \ for all \ \ } x \text{  and  } y,
\enas
for $\alpha \in (0,2/c)$ and $0 \le s \le \alpha$, we have
\beas
e^{sY^*}-e^{sY} \le \frac{|s(Y^*-Y)|(e^{sY^*}+e^{sY})}{2}.
\enas
Taking expectation and simplifying, we obtain
\beas
m^*(s) \le \left(\frac{1+cs/2}{1-cs/2}\right)m(s).
\enas
Using \eqref{m'mineq} and proceeding as in the first bound yield
\beas
m(\alpha) \le \exp\left(\frac{5(\sigma^2+B_2)\alpha^2}{6(1-c\alpha/2)} + B_1\alpha \right) \text{ \ \ for all \ \ } \alpha \in (0,2/c).
\enas
We note that the claim for $t=0$ is clear as the right hand side of \eqref{mainbound2} is one. Again applying Markov's inequality for $t>0$, we have
\beas
P(Y \ge t) \le \exp\left(-\alpha (t-B_1) + \frac{5(\sigma^2+B_2)\alpha^2}{6(1-c\alpha/2)}  \right)
\enas
Letting $\alpha = 2t/(10(\sigma^2+B_2)/2+ct)$ that lies in $(0,2/c)$, \eqref{mainbound2} follows.

(b) For any $s \in [0,\alpha)$ such that $m(\alpha)$ exists, using \eqref{m'bound}, \eqref{mrelate} and that $B_1=\E|R|/\lambda$ if $|R|$ and and $e^{sY}$ are negatively correlated for $0<s<\alpha$  and $B_1 = \esssup|\E[R|Y]|/\lambda$ otherwise and $B_2 = |\E(YR)|/\lambda$, we have
\beas
m'(s) &=& \sigma^2 s m^*(s) - \frac{\E YR}{\lambda} s m^*(s) +  \frac{\E [ R e^{sY} ]}{\lambda} \\
      &\le& \sigma^2 s m^*(s) - \frac{\E YR}{\lambda} s m^*(s) +  \frac{\E [e^{sY}|\E[R|Y]|]}{\lambda} \\
      &\le& (\sigma^2+B_2)se^{cs}m(s) + B_1m(s)
\enas
and thus
\beas
(\log m(s))' \le (\sigma^2+B_2)se^{cs} + B_1.
\enas

Integrating over $s \in [0,\alpha]$ we obtain
\beas
\log m(\alpha) \le \frac{\sigma^2+B_2}{c^2}(e^{c\alpha}(c\alpha-1)+1) + B_1 \alpha,
\enas
and hence
\beas
m(\alpha) \le \exp\left(\frac{\sigma^2+B_2}{c^2}(e^{c\alpha}(c\alpha-1)+1) + B_1 \alpha\right).
\enas

Again, by Markov's inequality,
\beas
P(Y \ge t) \le \exp\left(-\alpha(t-B_1)+\frac{\sigma^2+B_2}{c^2}(e^{c\alpha}(c\alpha-1)+1)\right)
\enas
Now, for $t>e+B_1$, choosing $\alpha = (\log (t-B_1) - \log \log (t-B_1))/c$, we have
\beas
P(Y \ge t) &\le& \exp\bigg(-\frac{t-B_1}{c}(\log (t-B_1) - \log \log (t-B_1)) \\
  &&+ \frac{\sigma^2+B_2}{c^2}\left(\frac{t-B_1}{\log (t-B_1) }(\log (t-B_1)  - \log \log (t-B_1) -1)+1\right)\bigg) \\
&\le& \exp\left(-\frac{t-B_1}{c}\left(\log (t-B_1)  - \log \log (t-B_1) -\frac{\sigma^2+B_2}{c}\right)\right).
\enas
Using that $(\log x)/2 \ge \log \log x$ for all $x>1$ yields the second inequality in \eqref{mainbound3}.

\bbox

\end{document}